\documentclass[a4paper, 12 pt]{article}

\usepackage[T2A]{fontenc}
\usepackage[cp1251]{inputenc}
\usepackage[english]{babel}
\usepackage[tbtags]{amsmath}
\usepackage{amsfonts,amssymb}

\usepackage{mathrsfs}

\usepackage{geometry} % Меняем поля страницы
\begin{document}
\title{$B$-expansion of  pseudo-involution in the Riordan group}
 \author{E. Burlachenko}
 \date{}

 \maketitle
\begin{abstract}
Each numerical sequence $\left( {{b}_{0}},{{b}_{1}},{{b}_{2}},... \right)$ with the generating function $B\left( x \right)$ defines the pseudo-involution in the Riordann group $\left( 1,xg\left( x \right) \right)$ such that $g\left( x \right)=1+xg\left( x \right)B\left( {{x}^{2}}g\left( x \right) \right)$. In the present paper we realize a simple idea: express the coefficients of the series ${{g}^{m}}\left( x \right)$ in terms of the coefficients of the series $B\left( x \right)$. Obtained  expansion  has a bright combinatorial character, sheds light on the connection of the pseudo-involution in the Riordann group with the generalized binomial series, and is also useful for finding  the series $g\left( x \right)$ by the given series $B\left( x \right)$. We compare this expansion with the similar expansion for the sequence $\left( 1,{{a}_{1}},{{a}_{2}},... \right)$  with the generating function $A\left( x \right)$ such that $g\left( x \right)=A\left( xg\left( x \right) \right)$.
\end{abstract}
\section{Introduction}

Transformations, corresponding to multiplication and composition of  series, play the main role in the space of formal power series over the field of real or complex  numbers. Multiplication is geven by the matrix$\left( f\left( x \right),x \right)$ $n$th column of which, $n=0,\text{ }1,\text{ }2,\text{ }...$ ,  has the generating function $f\left( x \right){{x}^{n}}$; composition is given by the matrix $\left( 1,g\left( x \right) \right)$ $n$th column of which  has the generating function ${{g}^{n}}\left( x \right)$, ${{g}_{0}}=0$:
$$\left( f\left( x \right),x \right)a\left( x \right)=f\left( x \right)a\left( x \right), \qquad\left( 1,g\left( x \right) \right)a\left( x \right)=a\left( g\left( x \right) \right).$$
Matrix
$$\left( f\left( x \right),x \right)\left( 1,g\left( x \right) \right)=\left( f\left( x \right),g\left( x \right) \right)$$
is called Riordan array [1] – [5]; $n$th column of Riordan array has the generating function $f\left( x \right){{g}^{n}}\left( x \right)$. Thus
$$\left( f\left( x \right),g\left( x \right) \right)b\left( x \right){{a}^{n}}\left( x \right)=f\left( x \right)b\left( g\left( x \right) \right){{a}^{n}}\left( g\left( x \right) \right),$$
$$\left( f\left( x \right),g\left( x \right) \right)\left( b\left( x \right),a\left( x \right) \right)=\left( f\left( x \right)b\left( g\left( x \right) \right),a\left( g\left( x \right) \right) \right).$$
Matrices $\left( f\left( x \right),g\left( x \right) \right)$, ${{f}_{0}}\ne 0$, ${{g}_{1}}\ne 0$, or in a more convenient  notation $\left( f\left( x \right),xg\left( x \right) \right)$, ${{f}_{0}}\ne 0$, ${{g}_{0}}\ne 0$, form a group, called the Riordan group. Elements of the matrix $\left( f\left( x \right),xg\left( x \right) \right)$ will be denoted ${{d}_{n,m}}$. For each matrix of the Riordan group there exists numerical sequence $A=\left( {{a}_{0}},{{a}_{1}},{{a}_{2}},... \right)$, called$A$-sequence, such that
$${{d}_{n+1,m+1}}=\sum\limits_{i=0}^{\infty }{{{a}_{i}}{{d}_{n,m+i}}}.$$
Let $A\left( x \right)$ is the generating function of the $A$-sequence. Then
$$f\left( x \right){{g}^{m+1}}\left( x \right)=f\left( x \right){{g}^{m}}\left( x \right)A\left( xg\left( x \right) \right),$$
$$g\left( x \right)=A\left( xg\left( x \right) \right),  \qquad{{\left( 1,xg\left( x \right) \right)}^{-1}}=\left( 1,x{{A}^{-1}}\left( x \right) \right).$$
For example,
$$A\left( x \right)=1+ax+b{{x}^{2}},  \qquad g\left( x \right)=1+axg\left( x \right)+b{{x}^{2}}{{g}^{2}}\left( x \right)=$$
$$=\frac{1-ax-\sqrt{{{\left( 1-ax \right)}^{2}}-4b{{x}^{2}}}}{2b{{x}^{2}}}.$$

Riordan array inverse to itself is called the involution in the Riordan group. If the matrix $\left( f\left( x \right),xg\left( x \right) \right)$ is an involution (in this case ${{f}_{0}},{{g}_{0}}=\pm 1$), then the matrix $\left( 1,xg\left( x \right) \right)$ is also an involution. The case ${{f}_{0}}=-1$ can be considered as the product of two involutions:
$$\left( -f\left( x \right),xg\left( x \right) \right)=\left( -1,x \right)\left( f\left( x \right),xg\left( x \right) \right), \qquad{{f}_{0}}=1.$$
Series $f\left( x \right)$, represented in the form
$$f\left( x \right)=c\left( x \right)+\sqrt{{{c}^{2}}\left( x \right)+1},  \qquad c\left( x \right)=\frac{f\left( x \right)-{{f}^{-1}}\left( x \right)}{2},$$
satisfies the condition 
$$c\left( xg\left( x \right) \right)=-c\left( x \right),  \qquad f\left( xg\left( x \right) \right)={{f}^{-1}}\left( x \right).$$
Any involution can be represented in the form $RM$, where $R$ is a Riordan array,
$$M=\left( 1,-x \right)=\left( \begin{matrix}
   1 & 0 & 0 & 0 & \cdots   \\
   0 & -1 & 0 & 0 & \cdots   \\
   0 & 0 & 1 & 0 & \cdots   \\
   0 & 0 & 0 & -1 & \cdots   \\
   \vdots  & \vdots  & \vdots  & \vdots  & \ddots   \\
\end{matrix} \right).$$
Matrix $R=\left( f\left( x \right),xg\left( x \right) \right)$,
$${{\left( f\left( x \right),xg\left( x \right) \right)}^{-1}}=M\left( f\left( x \right),xg\left( x \right) \right)M=\left( f\left( -x \right),xg\left( -x \right) \right),$$
is called the pseudo-involution in the Riordan group. An example of a pseudo-involution is the Pascal matrix:
$$P=\left( \frac{1}{1-x},\frac{x}{1-x} \right)=\left( \begin{matrix}
   1 & 0 & 0 & 0 & \cdots   \\
   1 & 1 & 0 & 0 & \cdots   \\
   1 & 2 & 1 & 0 & \cdots   \\
   1 & 3 & 3 & 1 & \cdots   \\
   \vdots  & \vdots  & \vdots  & \vdots  & \ddots   \\
\end{matrix} \right), \qquad{{P}^{-1}}=\left( \frac{1}{1+x},\frac{x}{1+x} \right).$$
For each pseudo-involution in the Riordan group (except matrices $M$, $-M$, which are simultaneously involutions and pseudo-involutions) there exists numerical sequence $B=\left( {{b}_{0}},{{b}_{1}},{{b}_{2}},... \right)$, called $B$-sequence [4], [5] (in [4] this sequence is called $\Delta $-sequence), such that
$${{d}_{n+1,m}}={{d}_{n,m-1}}+\sum\limits_{i=0}^{\infty }{{{b}_{i}}{{d}_{n-i,m+i}}}.$$
Let $B\left( x \right)$ is the generating function of the $B$-sequence of the matrix $\left( f\left( x \right),xg\left( x \right) \right)$. Then
$$f\left( x \right){{g}^{m}}\left( x \right)=f\left( x \right){{g}^{m-1}}\left( x \right)+xf\left( x \right){{g}^{m}}\left( x \right)B\left( {{x}^{2}}g\left( x \right) \right),$$
$$g\left( x \right)=1+xg\left( x \right)B\left( {{x}^{2}}g\left( x \right) \right).$$
For example,
$$B\left( x \right)=a+bx,  \qquad g\left( x \right)=1+axg\left( x \right)+b{{x}^{3}}{{g}^{2}}\left( x \right)=$$
$$=\frac{1-ax-\sqrt{{{\left( 1-ax \right)}^{2}}-4b{{x}^{3}}}}{2b{{x}^{3}}}.$$

In Section 2, for clarity which will be needed in the future, we associate the $B$-sequence of the matrix $\left( 1,xg\left( x \right) \right)$ with the $A$-sequence of the matrix $\left( 1,x\sqrt{g\left( x \right)} \right)$. In Section 3 on  basis of the identity
$${{g}^{m}}\left( x \right)={{g}^{m-1}}\left( x \right)+x{{g}^{m}}\left( x \right)B\left( {{x}^{2}}g\left( x \right) \right)$$
we express the coefficients of the series ${{g}^{m}}\left( x \right)$, $\left[ {{x}^{n}} \right]{{g}^{m}}\left( x \right)=g_{n}^{\left( m \right)}$, in terms of the coefficients of the series $B\left( x \right)$, namely
$$g_{n}^{\left( m \right)}=\sum\limits_{n}^{{}}{\frac{m{{\left( m+k \right)}_{q}}}{\left( m+k \right){{m}_{0}}!{{m}_{1}}!...{{m}_{p}}!}}b_{0}^{{{m}_{0}}}b_{1}^{{{m}_{1}}}...b_{p}^{{{m}_{p}}},$$ 
$$p=\left\lfloor \frac{n-1}{2} \right\rfloor ,  \qquad k=\sum\limits_{i=0}^{p}{{{m}_{i}}\left( i+1 \right)}, \qquad q=\sum\limits_{i=0}^{p}{{{m}_{i}}},$$
$${{\left( m+k \right)}_{q}}=\left( m+k \right)\left( m+k-1 \right)...\left( m+k-q+1 \right),$$
where the summation is over all monomials $b_{0}^{{{m}_{0}}}b_{1}^{{{m}_{1}}}...b_{p}^{{{m}_{p}}}$ for which $n=\sum\nolimits_{i=0}^{p}{{{m}_{i}}\left( 2i+1 \right)}$. In Section 4 we compare the obtained  expansion with expansions of the “binomial” and “generalized binomial” type, such as
 $$g_{n}^{\left( m \right)}=\sum\limits_{n}{\frac{{{\left( m \right)}_{q}}}{{{m}_{1}}!{{m}_{2}}...{{m}_{n}}!}}g_{1}^{{{m}_{1}}}g_{2}^{{{m}_{2}}}...g_{n}^{{{m}_{n}}}=\sum\limits_{n}{\frac{{{m}^{q}}}{{{m}_{1}}!{{m}_{2}}!\text{ }...\text{ }{{m}_{n}}!}}\text{ }l_{1}^{{{m}_{1}}}l_{2}^{{{m}_{2}}}...\text{ }l_{n}^{{{m}_{n}}}=$$
$$=\sum\limits_{n}{\frac{m{{\left( m+n \right)}_{q}}}{\left( m+n \right){{m}_{1}}!{{m}_{2}}...{{m}_{n}}!}}a_{1}^{{{m}_{1}}}a_{2}^{{{m}_{2}}}...a_{n}^{{{m}_{n}}},$$
 $${{l}_{i}}=\left[ {{x}^{i}} \right]\ln g\left( x \right), \qquad{{a}_{i}}=\left[ {{x}^{i}} \right]A\left( x \right),  \qquad n=\sum\limits_{i=1}^{n}{{{m}_{i}}i}, \qquad q=\sum\limits_{i=1}^{n}{{{m}_{i}}},$$
and show that it is also an expansion of this type.
\section{Some examples}
{\bfseries Remark 1.} If the matrices $\left( 1,x{{a}^{-1}}\left( x \right) \right)$, $\left( 1,xb\left( x \right) \right)$  are mutually inverse, then
$$\left( 1,xb\left( x \right) \right)a\left( x \right)=b\left( x \right),  \qquad\left( 1,xb\left( x \right) \right)\left( 1,xa\left( x \right) \right)=\left( 1,x{{b}^{2}}\left( x \right) \right).$$
Let ${{\left( 1,xa\left( x \right) \right)}^{-1}}=\left( 1,x{{c}^{-1}}\left( x \right) \right)$. Then
$$\left( 1,x{{c}^{-1}}\left( x \right) \right){{a}^{-1}}\left( x \right)={{c}^{-1}}\left( x \right),  \qquad\left( 1,x{{c}^{-1}}\left( x \right) \right)\left( 1,x{{a}^{-1}}\left( x \right) \right)=\left( 1,x{{c}^{-2}}\left( x \right) \right),$$
$${{\left( 1,x{{b}^{2}}\left( x \right) \right)}^{-1}}=\left( 1,x{{c}^{-2}}\left( x \right) \right).$$
{\bfseries Theorem 1.} \emph{If the matrix $\left( 1,xg\left( x \right) \right)$, $g\left( x \right)\ne -1$,  is a pseudo-involution, i.e. 
$${{\left( 1,xg\left( x \right) \right)}^{-1}}=\left( 1,xg\left( -x \right) \right)=M\left( 1,xg\left( x \right) \right)M,$$
then it can be represented in the form
$$\left( 1,xg\left( x \right) \right)=\left( 1,x\sqrt{g\left( x \right)} \right)\left( 1,xh\left( x \right) \right),$$
where  
$$h\left( -x \right)={{h}^{-1}}\left( x \right),   \qquad h\left( x \right)=s\left( x \right)+\sqrt{{{s}^{2}}\left( x \right)+1}, \qquad{{s}_{2n}}=0.$$ }
Proof follows from Remark 1.\\
{\bfseries Example 1.}
$$\left( 1,\frac{x}{1-2\varphi x} \right)=\left( 1,\frac{x}{\sqrt{1-2\varphi x}} \right)\left( 1,x\left( \varphi x+\sqrt{{{\varphi }^{2}}{{x}^{2}}+1} \right) \right).$$
{\bfseries Example 2.}
$$\left( 1,x\sum\limits_{n=0}^{\infty }{\frac{2{{\left( 2+n \right)}^{n-1}}}{n!}{{\varphi }^{n}}{{x}^{n}}} \right)=\left( 1,x\sum\limits_{n=0}^{\infty }{\frac{{{\left( 1+n \right)}^{n-1}}}{n!}{{\varphi }^{n}}{{x}^{n}}} \right)\left( 1,x{{e}^{\varphi x}} \right),$$
where
$$x\sum\limits_{n=0}^{\infty }{\frac{{{\left( 1+n \right)}^{n-1}}}{n!}{{\varphi }^{n}}{{x}^{n}}=\ln \left( \sum\limits_{n=0}^{\infty }{\frac{{{\left( 1+\varphi n \right)}^{n-1}}}{n!}{{x}^{n}}} \right)}=x{{\left( \sum\limits_{n=0}^{\infty }{\frac{{{\left( 1+\varphi n \right)}^{n-1}}}{n!}{{x}^{n}}} \right)}^{\varphi }},$$
$$\sum\limits_{n=0}^{\infty }{\frac{2{{\left( 2+n \right)}^{n-1}}}{n!}{{\varphi }^{n}}{{x}^{n}}={{\left( \sum\limits_{n=0}^{\infty }{\frac{{{\left( 1+\varphi n \right)}^{n-1}}}{n!}{{x}^{n}}} \right)}^{2\varphi }}}.$$
{\bfseries Example 3.}
$$\left( 1,\frac{1-4\varphi x+{{\varphi }^{2}}{{x}^{2}}-\sqrt{{{\left( 1-4\varphi x+{{\varphi }^{2}}{{x}^{2}} \right)}^{2}}-4{{\varphi }^{2}}{{x}^{2}}}}{2{{\varphi }^{2}}x} \right)=$$
$$=\left( 1,\frac{1-\varphi x-\sqrt{{{\left( 1-\varphi x \right)}^{2}}-4\varphi x}}{2\varphi } \right)\left( 1,x\frac{1+\varphi x}{1-\varphi x} \right),$$
$$\frac{1+\varphi x}{1-\varphi x}=\frac{2\varphi x}{1-{{\varphi }^{2}}{{x}^{2}}}+\sqrt{{{\left( \frac{2\varphi x}{1-{{\varphi }^{2}}{{x}^{2}}} \right)}^{2}}+1}.$$
{\bfseries Theorem 2.} \emph{If $B\left( x \right)$ is the generating function of the $B$-sequence of the matrix $\left( 1,xg\left( x \right) \right)$, then
$$xB\left( {{x}^{2}} \right)=2s\left( x \right).$$}
{\bfseries Proof.} Since ${{h}^{2}}\left( x \right)=1+2s\left( x \right)h\left( x \right)$, then
$$g\left( x \right)=\left( 1,x\sqrt{g\left( x \right)} \right)\left( 1+2s\left( x \right)h\left( x \right) \right)=1+xg\left( x \right)\tilde{s}\left( x\sqrt{g\left( x \right)} \right)=$$
$$=1+xg\left( x \right)B\left( {{x}^{2}}g\left( x \right) \right),  \qquad\tilde{s}\left( x \right)=\frac{2s\left( x \right)}{x}.$$
{\bfseries Example 4.}
 Paper [5] contains the interesting fact that if
$$g\left( x \right)=\sum\limits_{n=0}^{\infty }{\frac{2m+1}{2m+1+\left( m+1 \right)n}}\left( \begin{matrix}
   2m+1+\left( m+1 \right)n  \\
   n  \\
\end{matrix} \right){{x}^{n}},$$
then $B$-sequence of the matrix $\left( 1,xg\left( x \right) \right)$ coincides with the $m$th row of the matrix 
$$\left( \frac{1+x}{{{\left( 1-x \right)}^{2}}},\frac{x}{{{\left( 1-x \right)}^{2}}} \right)=\left( \begin{matrix}
   1 & 0 & 0 & 0 & \cdots   \\
   3 & 1 & 0 & 0 & \cdots   \\
   5 & 5 & 1 & 0 & \cdots   \\
   7 & 14 & 7 & 1 & \cdots   \\
   \vdots  & \vdots  & \vdots  & \vdots  & \ddots   \\
\end{matrix} \right).$$
This is consequence of the fact that in this case
$$h\left( x \right)={{\left( \frac{x+\sqrt{{{x}^{2}}+4}}{2} \right)}^{2m+1}},$$
$${{\left( \frac{x+\sqrt{{{x}^{2}}+4}}{2} \right)}^{n}}=\frac{{{c}_{n}}\left( x \right)+{{s}_{n-1}}\left( x \right)\sqrt{{{x}^{2}}+4}}{2},$$
$${{s}_{2m}}\left( x \right)\sqrt{{{x}^{2}}+4}=\sqrt{c_{2m+1}^{2}\left( x \right)+4},  \qquad{{c}_{2m}}\left( x \right)=\sqrt{s_{2m-1}^{2}\left( x \right)\left( {{x}^{2}}+4 \right)+4},$$
where polynomial ${{c}_{n}}\left( x \right)$ corresponds to the $n$th row of the matrix
$$\left( \frac{1+{{x}^{2}}}{1-{{x}^{2}}},\frac{x}{1-{{x}^{2}}} \right)=\left( \begin{matrix}
   1 & 0 & 0 & 0 & 0 & 0 & \cdots   \\
   0 & 1 & 0 & 0 & 0 & 0 & \cdots   \\
   2 & 0 & 1 & 0 & 0 & 0 & \cdots   \\
   0 & 3 & 0 & 1 & 0 & 0 & \cdots   \\
   2 & 0 & 4 & 0 & 1 & 0 & \cdots   \\
   0 & 5 & 0 & 5 & 0 & 1 & \cdots   \\
   \vdots  & \vdots  & \vdots  & \vdots  & \vdots  & \vdots  & \ddots   \\
\end{matrix} \right),$$
polynomial ${{s}_{n}}\left( x \right)$ corresponds to the $n$th row of the matrix
$$\left( \frac{1}{1-{{x}^{2}}},\frac{x}{1-{{x}^{2}}} \right)=\left( \begin{matrix}
   1 & 0 & 0 & 0 & 0 & 0 & \cdots   \\
   0 & 1 & 0 & 0 & 0 & 0 & \cdots   \\
   1 & 0 & 1 & 0 & 0 & 0 & \cdots   \\
   0 & 2 & 0 & 1 & 0 & 0 & \cdots   \\
   1 & 0 & 3 & 0 & 1 & 0 & \cdots   \\
   0 & 3 & 0 & 4 & 0 & 1 & \cdots   \\
   \vdots  & \vdots  & \vdots  & \vdots  & \vdots  & \vdots  & \ddots   \\
\end{matrix} \right).$$
\section{$B$-expansion}
Denote  $\left[ {{x}^{n}} \right]{{g}^{m}}\left( x \right)=g_{n}^{\left( m \right)}$, $g_{n}^{\left( 1 \right)}={{g}_{n}}$ . Since
$${{g}^{m}}\left( x \right)={{g}^{m-1}}\left( x \right)+x{{g}^{m}}\left( x \right)B\left( {{x}^{2}}g\left( x \right) \right),$$
where
$${{g}^{m}}\left( x \right)B\left( {{x}^{2}}g\left( x \right) \right)=\left( {{g}^{m}}\left( x \right),{{x}^{2}}g\left( x \right) \right)B\left( x \right),$$
then
$$g_{n}^{\left( m \right)}={{b}_{0}}g_{n-1}^{\left( m \right)}+{{b}_{1}}g_{n-3}^{\left( m+1 \right)}+{{b}_{2}}g_{n-5}^{\left( m+2 \right)}+...+{{b}_{p}}g_{n-1-2p}^{\left( m+p \right)}+g_{n}^{\left( m-1 \right)}, \qquad p=\left\lfloor \frac{n-1}{2} \right\rfloor ,$$
or
$$g_{n}^{\left( m \right)}={{b}_{0}}g_{n-1}^{\left( m \right)}+{{b}_{1}}g_{n-3}^{\left( m+1 \right)}+{{b}_{2}}g_{n-5}^{\left( m+2 \right)}+...+{{b}_{p}}g_{n-1-2p}^{\left( m+p \right)}+$$
$$+{{b}_{0}}g_{n-1}^{\left( m-1 \right)}+{{b}_{1}}g_{n-3}^{\left( m \right)}+{{b}_{2}}g_{n-5}^{\left( m+1 \right)}+...+{{b}_{p}}g_{n-1-2p}^{\left( m+p-1 \right)}+$$
$$+{{b}_{0}}g_{n-1}^{\left( m-2 \right)}+{{b}_{1}}g_{n-3}^{\left( m-1 \right)}+{{b}_{2}}g_{n-5}^{\left( m \right)}+...+{{b}_{p}}g_{n-1-2p}^{\left( m+p-2 \right)}+$$
$$…$$
$$+{{b}_{0}}g_{n-1}^{\left( 1 \right)}+{{b}_{1}}g_{n-3}^{\left( 2 \right)}+{{b}_{2}}g_{n-5}^{\left( 3 \right)}+...+{{b}_{p}}g_{n-1-2p}^{\left( p+1 \right)};$$
$$g_{n}^{\left( m \right)}={{b}_{0}}\sum\limits_{i=1}^{m}{g_{n-1}^{\left( i \right)}}+{{b}_{1}}\sum\limits_{i=2}^{m+1}{g_{n-3}^{\left( i \right)}}+{{b}_{2}}\sum\limits_{i=3}^{m+2}{g_{n-5}^{\left( i \right)}+...+{{b}_{p}}}\sum\limits_{i=p+1}^{m+p}{g_{n-1-2p}^{\left( i \right)}},\eqno (1)$$
Using recursion, we find
$$g_{0}^{\left( m \right)}=1,  \qquad g_{1}^{\left( m \right)}=m{{b}_{0}},  \qquad g_{2}^{\left( m \right)}=\left( \begin{matrix}
   m+1  \\
   2  \\
\end{matrix} \right)b_{0}^{2},$$
 $$g_{3}^{\left( m \right)}=\left( \begin{matrix}
   m+2  \\
   3  \\
\end{matrix} \right)b_{0}^{3}+m{{b}_{1}},$$
$$g_{4}^{\left( m \right)}=\left( \begin{matrix}
   m+3  \\
   4  \\
\end{matrix} \right)b_{0}^{4}+m\left( \begin{matrix}
   m+2  \\
   1  \\
\end{matrix} \right){{b}_{0}}{{b}_{1}},$$
$$g_{5}^{\left( m \right)}=\left( \begin{matrix}
   m+4  \\
   5  \\
\end{matrix} \right)b_{0}^{5}+m\left( \begin{matrix}
   m+3  \\
   2  \\
\end{matrix} \right)b_{0}^{2}{{b}_{1}}+m{{b}_{2}},$$
$$g_{6}^{\left( m \right)}=\left( \begin{matrix}
   m+5  \\
   6  \\
\end{matrix} \right)b_{0}^{6}+m\left( \begin{matrix}
   m+4  \\
   3  \\
\end{matrix} \right)b_{0}^{3}{{b}_{1}}+m\left( \begin{matrix}
   m+3  \\
   1  \\
\end{matrix} \right){{b}_{0}}{{b}_{2}}+\frac{m}{m+2}\left( \begin{matrix}
   m+3  \\
   2  \\
\end{matrix} \right)b_{1}^{2}.$$
Coefficient of the monomial $b_{0}^{{{m}_{0}}}b_{1}^{{{m}_{1}}}...b_{p}^{{{m}_{p}}}$ in the expansion of the coefficient $g_{n}^{\left( m \right)}$ will be denoted $\left( m|b_{0}^{{{m}_{0}}}b_{1}^{{{m}_{1}}}...b_{p}^{{{m}_{p}}} \right)$.\\
{\bfseries Theorem 3.}\emph{$$g_{n}^{\left( m \right)}=\sum{\left( m|b_{0}^{{{m}_{0}}}b_{1}^{{{m}_{1}}}...b_{p}^{{{m}_{p}}} \right)}b_{0}^{{{m}_{0}}}b_{1}^{{{m}_{1}}}...b_{p}^{{{m}_{p}}},$$
where expression $b_{0}^{{{m}_{0}}}b_{1}^{{{m}_{1}}}...b_{p}^{{{m}_{p}}}$ corresponds to the partition $n=\sum\nolimits_{i=0}^{p}{{{m}_{i}}}\left( 2i+1 \right)$ and summation is done over all partitions of the number $n$ into odd parts.}\\
{\bfseries Proof.} Let the theorem is true for ${{g}_{n}}$:
$${{g}_{n}}=\sum{\left( 1|b_{0}^{{{m}_{0}}}b_{1}^{{{m}_{1}}}...b_{p}^{{{m}_{p}}} \right)}b_{0}^{{{m}_{0}}}b_{1}^{{{m}_{1}}}...b_{p}^{{{m}_{p}}}.$$
Then
$$g_{n}^{\left( 2 \right)}=\sum\limits_{n=i+j}{{{g}_{i}}{{g}_{j}}}=\sum{\left( 2|b_{0}^{{{m}_{0}}}b_{1}^{{{m}_{1}}}...b_{p}^{{{m}_{p}}} \right)}b_{0}^{{{m}_{0}}}b_{1}^{{{m}_{1}}}...b_{p}^{{{m}_{p}}},$$
$$g_{n}^{\left( m \right)}=\sum\limits_{n=i+j}{{{g}_{i}}g_{j}^{\left( m-1 \right)}}=\sum{\left( m|b_{0}^{{{m}_{0}}}b_{1}^{{{m}_{1}}}...b_{p}^{{{m}_{p}}} \right)}b_{0}^{{{m}_{0}}}b_{1}^{{{m}_{1}}}...b_{p}^{{{m}_{p}}}.$$
I.e. set of monomials in the expansion of the coefficient $g_{n}^{\left( m \right)}$ does not depend on $m$. Let the theorem is true for all ${{g}_{i}}$, $i<n$. Then it is also true for ${{g}_{n}}$,
$${{g}_{n}}={{b}_{0}}g_{n-1}^{\left( 1 \right)}+{{b}_{1}}g_{n-3}^{\left( 2 \right)}+{{b}_{2}}g_{n-5}^{\left( 3 \right)}+...+{{b}_{p}}g_{n-1-2p}^{\left( p+1 \right)},$$
since monomial, corresponding to the partition $n=\sum\nolimits_{i=0}^{p}{{{m}_{i}}}\left( 2i+1 \right)$,  is contained in the summand ${{b}_{i}}g_{n-1-2i}^{\left( i+1 \right)}$, if ${{m}_{i}}\ne 0$. Thus, it is sufficient that the theorem was true for ${{g}_{1}}$.\\
{\bfseries Theorem 4.}\emph{$$\left( m|b_{0}^{{{m}_{0}}}b_{1}^{{{m}_{1}}}...b_{p}^{{{m}_{p}}} \right)=\sum\limits_{i=1}^{m}{\left( i|b_{0}^{{{m}_{0}}-1}b_{1}^{{{m}_{1}}}...b_{p}^{{{m}_{p}}} \right)}+\sum\limits_{i=2}^{m+1}{\left( i|b_{0}^{{{m}_{0}}}b_{1}^{{{m}_{1}}-1}...b_{p}^{{{m}_{p}}} \right)}+$$
$$+\sum\limits_{i=3}^{m+2}{\left( i|b_{0}^{{{m}_{0}}}b_{1}^{{{m}_{1}}}b_{2}^{{{m}_{2}}-1}...b_{p}^{{{m}_{p}}} \right)}+...+\sum\limits_{i=p+1}^{m+p}{\left( i|b_{0}^{{{m}_{0}}}b_{1}^{{{m}_{1}}}...b_{p}^{{{m}_{p}}-1} \right)},$$ 
where $\left( i|...b_{r}^{-1}... \right)=0$.}\\
{\bfseries Proof.} From Theorem 3 it follows that the monomial $b_{0}^{{{m}_{0}}}b_{1}^{{{m}_{1}}}b_{2}^{{{m}_{2}}}...b_{p}^{{{m}_{p}}}$ with the coefficient $\sum\nolimits_{i=r+1}^{m+r}{\left( i|b_{0}^{{{m}_{0}}}b_{1}^{{{m}_{1}}}...b_{r}^{{{m}_{r}}-1}...b_{p}^{{{m}_{p}}} \right)}$  is present in the summand 
$${{b}_{r}}\sum\limits_{i=r+1}^{m+r}{g_{n-1-2r}^{\left( i \right)}}$$
of the formula (1), if ${{m}_{r}}\ne 0$.

Coefficients $\left( m|b_{0}^{{{m}_{0}}}b_{1}^{{{m}_{1}}}...b_{p}^{{{m}_{p}}} \right)$ are closely related to the coefficients of the generalized binomial series
$${{\mathcal{B}}_{r}}{{\left( x \right)}^{m}}=\sum\limits_{n=0}^{\infty }{\frac{m}{m+rn}}\left( \begin{matrix}
   m+rn  \\
   n  \\
\end{matrix} \right){{x}^{n}}.$$
Consider the following generalization of the Pascal table. Elements of the table will be denoted ${{\left( m,n \right)}_{r}}$. Then ${{\left( m,0 \right)}_{r}}=1$; ${{\left( 0,n \right)}_{r}}=0$, $n>0$. Remaining elements will be found by the rule 
$${{\left( m,n \right)}_{r}}={{\left( m-1,n \right)}_{r}}+{{\left( m+r-1,n-1 \right)}_{r}}.$$
For example, $r=1$, $r=2$, $r=3$, $r=4$:
$$\left( \begin{matrix}
   1 & 0 & 0 & 0 & 0 & \cdots   \\
   1 & 1 & 1 & 1 & 1 & \cdots   \\
   1 & 2 & 3 & 4 & 5 & \cdots   \\
   1 & 3 & 6 & 10 & 15 & \cdots   \\
   1 & 4 & 10 & 20 & 35 & \cdots   \\
   \vdots  & \vdots  & \vdots  & \vdots  & \vdots  & \ddots   \\
\end{matrix} \right),  \qquad\left( \begin{matrix}
   1 & 0 & 0 & 0 & 0 & \cdots   \\
   1 & 1 & 2 & 5 & 14 & \cdots   \\
   1 & 2 & 5 & 14 & 42 & \cdots   \\
   1 & 3 & 9 & 28 & 92 & \cdots   \\
   1 & 4 & 14 & 48 & 165 & \cdots   \\
   \vdots  & \vdots  & \vdots  & \vdots  & \vdots  & \ddots   \\
\end{matrix} \right),$$
$$\left( \begin{matrix}
   1 & 0 & 0 & 0 & 0 & \cdots   \\
   1 & 1 & 3 & 12 & 55 & \cdots   \\
   1 & 2 & 7 & 30 & 143 & \cdots   \\
   1 & 3 & 12 & 55 & 273 & \cdots   \\
   1 & 4 & 18 & 88 & 455 & \cdots   \\
   \vdots  & \vdots  & \vdots  & \vdots  & \vdots  & \ddots   \\
\end{matrix} \right), \qquad\left( \begin{matrix}
   1 & 0 & 0 & 0 & 0 & \cdots   \\
   1 & 1 & 4 & 22 & 140 & \cdots   \\
   1 & 2 & 9 & 52 & 340 & \cdots   \\
   1 & 3 & 15 & 91 & 612 & \cdots   \\
   1 & 4 & 22 & 140 & 967 & \cdots   \\
   \vdots  & \vdots  & \vdots  & \vdots  & \vdots  & \ddots   \\
\end{matrix} \right).$$
Then 
$${{\left( m,n \right)}_{r}}=\sum\limits_{i=r}^{m+r-1}{{{\left( i,n-1 \right)}_{r}}},   \qquad{{\mathcal{B}}_{r}}{{\left( x \right)}^{m}}=\sum\limits_{n=0}^{\infty }{{{\left( m,n \right)}_{r}}}{{x}^{n}}.$$\\
{\bfseries Theorem 5.}
$$\left( m|b_{r}^{{{m}_{r}}} \right)=\left[ {{x}^{{{m}_{r}}}} \right]{{\mathcal{B}}_{r+1}}{{\left( x \right)}^{m}}=\frac{m}{m+r{{m}_{r}}}\left( \begin{matrix}
   m+r{{m}_{r}}+{{m}_{r}}-1  \\
   {{m}_{r}}  \\
\end{matrix} \right).$$
{\bfseries Proof.} According to the Theorem 4
$$\left( m|b_{r}^{{{m}_{r}}} \right)=\sum\limits_{i=r+1}^{m+r}{\left( i|b_{r}^{{{m}_{r}}-1} \right)},$$
where $\left( i|{{b}_{r}} \right)=i$.\\
{\bfseries Theorem 6.}
$$\left( m|b_{r}^{{{m}_{r}}}b_{s}^{{{m}_{s}}} \right)=\frac{m}{m+k-{{m}_{r}}-{{m}_{s}}}\left( \begin{matrix}
   m+k-1  \\
   {{m}_{r}}  \\
\end{matrix} \right)\left( \begin{matrix}
   m+k-1-{{m}_{r}}  \\
   {{m}_{s}}  \\
\end{matrix} \right)=$$
$$=\frac{m\left( m+k-1 \right)!}{{{m}_{r}}!{{m}_{s}}!\left( m+k-{{m}_{r}}-{{m}_{s}} \right)!},  \qquad k={{m}_{r}}\left( r+1 \right)+{{m}_{s}}\left( s+1 \right).$$
{\bfseries Proof.} By successively applying Theorem 4, we can expand the coefficients $\left( m|b_{r}^{{{m}_{r}}}b_{s}^{{{m}_{s}}} \right)$ into a sum of the coefficients of the form $\left( i|b_{j}^{{{m}_{j}}} \right)$ which satisfy Theorem 6. Therefore it suffices to show that Theorem 4 is compatible with Theorem 6:
$$\left( m|b_{r}^{{{m}_{r}}}b_{s}^{{{m}_{s}}} \right)=\sum\limits_{i=r+1}^{m+r}{\left( i|b_{r}^{{{m}_{r}}-1}b_{s}^{{{m}_{s}}} \right)}+\sum\limits_{i=s+1}^{m+s}{\left( i|b_{r}^{{{m}_{r}}}b_{s}^{{{m}_{s}}-1} \right)}=$$
$$=\sum\limits_{i=r+1}^{m+r}{\frac{i\left( i+k-1-r-1 \right)!}{\left( {{m}_{r}}-1 \right)!{{m}_{s}}!\left( i+k-{{m}_{r}}-{{m}_{s}}-r \right)!}+}$$
$$+\sum\limits_{i=s+1}^{m+s}{\frac{i\left( i+k-1-s-1 \right)!}{{{m}_{r}}!\left( {{m}_{s}}-1 \right)!\left( i+k-{{m}_{r}}-{{m}_{s}}-s \right)!}}=$$
$$=\sum\limits_{i=1}^{m}{\frac{\left( {{m}_{r}}\left( r+i \right)+{{m}_{s}}\left( s+i \right) \right)\left( i+k-2 \right)!}{{{m}_{r}}!{{m}_{s}}!\left( i+k-{{m}_{r}}-{{m}_{s}} \right)!}}=$$ 
$$=\frac{k!}{{{m}_{r}}!{{m}_{s}}!\left( 1+k-{{m}_{r}}-{{m}_{s}} \right)!}+\sum\limits_{i=2}^{m}{\frac{\left( {{m}_{r}}\left( r+i \right)+{{m}_{s}}\left( s+i \right) \right)\left( i+k-2 \right)!}{{{m}_{r}}!{{m}_{s}}!\left( i+k-{{m}_{r}}-{{m}_{s}} \right)!}}=$$
$$=\frac{2\left( 1+k \right)!}{{{m}_{r}}!{{m}_{s}}!\left( 2+k-{{m}_{r}}-{{m}_{s}} \right)!}+\sum\limits_{i=3}^{m}{\frac{\left( {{m}_{r}}\left( r+i \right)+{{m}_{s}}\left( s+i \right) \right)\left( i+k-2 \right)!}{{{m}_{r}}!{{m}_{s}}!\left( i+k-{{m}_{r}}-{{m}_{s}} \right)!}}=$$

$$…$$

$$=\frac{\left( m-1 \right)\left( m+k-2 \right)!}{{{m}_{r}}!{{m}_{s}}!\left( m-1+k-{{m}_{r}}-{{m}_{s}} \right)!}+\frac{\left( {{m}_{r}}\left( r+m \right)+{{m}_{s}}\left( s+m \right) \right)\left( m+k-2 \right)!}{{{m}_{r}}!{{m}_{s}}!\left( m+k-{{m}_{r}}-{{m}_{s}} \right)!}=$$
$$=\frac{\left( \left( m-1 \right)\left( m+k \right)+k \right)\left( m+k-2 \right)!}{{{m}_{r}}!{{m}_{s}}!\left( m+k-{{m}_{r}}-{{m}_{s}} \right)!}=\frac{m\left( m+k-1 \right)!}{{{m}_{r}}!{{m}_{s}}!\left( m+k-{{m}_{r}}-{{m}_{s}} \right)!}.$$
Generalizing, we deduce
$$\left( m|b_{0}^{{{m}_{0}}}b_{1}^{{{m}_{1}}}...b_{p}^{{{m}_{p}}} \right)=\frac{m\left( m+k-1 \right)!}{{{m}_{0}}!{{m}_{1}}!...{{m}_{p}}!\left( m+k-{{m}_{0}}-{{m}_{1}}-...-{{m}_{p}} \right)!},$$
$$k=\sum\limits_{i=0}^{p}{{{m}_{i}}\left( i+1 \right)}.$$

Let the expression 
$$\sum\limits_{n}^{{}}{\left( m|b_{0}^{{{m}_{0}}}b_{1}^{{{m}_{1}}}...b_{p}^{{{m}_{p}}} \right)b_{0}^{{{m}_{0}}}b_{1}^{{{m}_{1}}}...b_{p}^{{{m}_{p}}}}$$
mean that the summation is over all monomials $b_{0}^{{{m}_{0}}}b_{1}^{{{m}_{1}}}...b_{p}^{{{m}_{p}}}$ for which $n=\sum\nolimits_{i=0}^{p}{{{m}_{i}}\left( 2i+1 \right)}$ (or by another rule for $n$, which is indicated separately). Then
$${{g}^{m}}\left( x \right)=1+\sum\limits_{n=1}^{\infty }{\sum\limits_{n}^{{}}{\left( m|b_{0}^{{{m}_{0}}}b_{1}^{{{m}_{1}}}...b_{p}^{{{m}_{p}}} \right)b_{0}^{{{m}_{0}}}b_{1}^{{{m}_{1}}}...b_{p}^{{{m}_{p}}}{{x}^{n}}}}.$$ 
Since
$$\left( m|b_{0}^{{{m}_{0}}}b_{1}^{{{m}_{1}}}b_{2}^{{{m}_{2}}}...b_{p}^{{{m}_{p}}} \right)=$$
$$=\frac{\left( m+k-1 \right)!m\left( m+k-{{m}_{0}}-1 \right)!}{{{m}_{0}}!\left( m+k-{{m}_{0}}-1 \right)!{{m}_{1}}!...{{m}_{p}}!\left( m+k-{{m}_{0}}-{{m}_{1}}-...-{{m}_{p}} \right)!}=$$
$$=\left( \begin{matrix}
   m+\left( k-{{m}_{0}} \right)+{{m}_{0}}-1  \\
   {{m}_{0}}  \\
\end{matrix} \right)\left( m|b_{1}^{{{m}_{1}}}b_{2}^{{{m}_{2}}}...b_{p}^{{{m}_{p}}} \right),$$
then the series ${{g}^{m}}\left( x \right)$ can also be represented in the form
$${{g}^{m}}\left( x \right)=\frac{1}{{{\left( 1-{{b}_{0}}x \right)}^{m}}}+\sum\limits_{n=1}^{\infty }{\sum\limits_{n}^{{}}{\left( m|b_{1}^{{{m}_{1}}}b_{2}^{{{m}_{2}}}...b_{p}^{{{m}_{p}}} \right)b_{1}^{{{m}_{1}}}b_{2}^{{{m}_{2}}}...b_{p}^{{{m}_{p}}}\frac{{{x}^{n}}}{{{\left( 1-{{b}_{0}}x \right)}^{m+k}}}}},$$
$$n=\sum\limits_{i=1}^{p}{{{m}_{i}}}\left( 2i+1 \right), \qquad k=\sum\limits_{i=1}^{p}{{{m}_{i}}}\left( i+1 \right).$$
{\bfseries Example 5.}
$$B\left( x \right)={{b}_{0}}+{{b}_{r}}{{x}^{r}},  \qquad{{g}^{m}}\left( x \right)=\frac{1}{{{\left( 1-{{b}_{0}}x \right)}^{m}}}+\sum\limits_{n=1}^{\infty }{\left( m|b_{r}^{n} \right)b_{r}^{n}\frac{{{x}^{n\left( 2r+1 \right)}}}{{{\left( 1-{{b}_{0}}x \right)}^{m+n\left( r+1 \right)}}}}=$$
$$=\left( \frac{1}{{{\left( 1-{{b}_{0}}x \right)}^{m}}},\frac{{{b}_{r}}{{x}^{2r+1}}}{{{\left( 1-{{b}_{0}}x \right)}^{r+1}}} \right){{\mathcal{B}}_{r+1}}{{\left( x \right)}^{m}}.$$
In particular,
$$\left( \frac{1}{{{\left( 1-{{b}_{0}}x \right)}^{m}}},\frac{{{b}_{1}}{{x}^{3}}}{{{\left( 1-{{b}_{0}}x \right)}^{2}}} \right){{\left( \frac{1-\sqrt{1-4x}}{2x} \right)}^{m}}
={{\left( \frac{1-{{b}_{0}}x-\sqrt{{{\left( 1-{{b}_{0}}x \right)}^{2}}-4{{b}_{1}}{{x}^{3}}}}{2{{b}_{1}}{{x}^{3}}} \right)}^{m}}.$$
\section{Expansions of  generalized binomial type}
Let $\left| {{e}^{x}} \right|$ is the diagonal matrix whose diagonal elements are equal to the coefficients of the series ${{e}^{x}}$: $\left| {{e}^{x}} \right|{{\left( 1-x \right)}^{-1}}={{e}^{x}}$. Polynomial, corresponding to the $n$th row of the matrix ${{\left| {{e}^{x}} \right|}^{-1}}\left( 1,\ln a\left( x \right) \right)\left| {{e}^{x}} \right|$, will be denoted ${{p}_{n}}\left( x \right)$ (sequence of such polynomials is called the binomial sequence). Then
$${{a}^{\varphi }}\left( x \right)=\sum\limits_{n=0}^{\infty }{\frac{{{p}_{n}}\left( \varphi  \right)}{n!}{{x}^{n}}}.$$
Polynomial, corresponding to the $n$th row of the matrix $\left( 1,f\left( x \right) \right)$, ${{f}_{0}}=0$, ${{f}_{1}}\ne 0$,    $n>0$, has the form
$$\sum\limits_{n}{\frac{q!{{x}^{q}}}{{{m}_{1}}!{{m}_{2}}!...{{m}_{n}}!}f_{1}^{{{m}_{1}}}f_{2}^{{{m}_{2}}}...f_{n}^{{{m}_{n}}}},  \qquad n=\sum\limits_{i=1}^{n}{{{m}_{i}}}i,  \qquad q=\sum\limits_{i=1}^{n}{{{m}_{i}}}.$$
Hence, if $g\left( x \right)=a\left( f\left( x \right) \right)$, then
$$g_{n}^{\left( m \right)}=\sum\limits_{n}{\frac{{{p}_{q}}\left( m \right)}{{{m}_{1}}!{{m}_{2}}!...{{m}_{n}}!}f_{1}^{{{m}_{1}}}f_{2}^{{{m}_{2}}}...f_{n}^{{{m}_{n}}}}.$$
Representation of the coefficients $g_{n}^{\left( m \right)}$ in this form will be called expansion of the binomial type, or the binomial expansion. For example, since
$${{g}^{m}}\left( x \right)=\left( 1,g\left( x \right)-1 \right){{\left( 1+x \right)}^{m}}=\left( 1,\ln g\left( x \right) \right){{e}^{xm}},$$
then
$$g_{n}^{\left( m \right)}=\sum\limits_{n}{\frac{{{\left( m \right)}_{q}}}{{{m}_{1}}!{{m}_{2}}...{{m}_{n}}!}}g_{1}^{{{m}_{1}}}g_{2}^{{{m}_{2}}}...g_{n}^{{{m}_{n}}}=\sum\limits_{n}{\frac{{{m}^{q}}}{{{m}_{1}}!{{m}_{2}}!\text{ }...\text{ }{{m}_{n}}!}}\text{ }l_{1}^{{{m}_{1}}}l_{2}^{{{m}_{2}}}...\text{ }l_{n}^{{{m}_{n}}},$$
$${{l}_{n}}=\left[ {{x}^{n}} \right]\ln g\left( x \right),  \qquad n=\sum\limits_{i=1}^{n}{{{m}_{i}}}i ,  \qquad q=\sum\limits_{i=1}^{n}{{{m}_{i}}}.$$
Polynomial, corresponding to the $n$th row of the matrix $\left( 1,\ln g\left( x \right) \right)\left| {{e}^{x}} \right|$ , will be denoted ${{l}_{n}}\left( x \right)$. Then
$${{l}_{n}}\left( x \right)=\sum\limits_{n}{\frac{{{p}_{q}}\left( x \right)}{{{m}_{1}}!{{m}_{2}}!...{{m}_{n}}!}f_{1}^{{{m}_{1}}}f_{2}^{{{m}_{2}}}...f_{n}^{{{m}_{n}}}}=\sum\limits_{n}{\frac{{{x}^{q}}}{{{m}_{1}}!{{m}_{2}}!\text{ }...\text{ }{{m}_{n}}!}}\text{ }l_{1}^{{{m}_{1}}}l_{2}^{{{m}_{2}}}...\text{ }l_{n}^{{{m}_{n}}}.$$
Polynomial, corresponding to the $n$th row of the matrix $\left( 1,\ln A\left( x \right) \right)\left| {{e}^{x}} \right|$, will be denoted ${{\tilde{l}}_{n}}\left( x \right)$. Since ${{\left( 1,xg\left( x \right) \right)}^{-1}}=\left( 1,x{{A}^{-1}}\left( x \right) \right)$, then by the Lagrange inversion theorem
$${{l}_{n}}\left( x \right)=x{{\left( x+n \right)}^{-1}}{{\tilde{l}}_{n}}\left( x+n \right).$$
Thus,
$${{\tilde{l}}_{n}}\left( x \right)=\sum\limits_{n}^{{}}{\frac{{{\left( x \right)}_{q}}}{{{m}_{1}}!{{m}_{2}}...{{m}_{n}}!}}a_{1}^{{{m}_{1}}}a_{2}^{{{m}_{2}}}...a_{n}^{{{m}_{n}}},$$
$${{l}_{n}}\left( x \right)=\sum\limits_{n}^{{}}{\frac{x{{\left( x+n \right)}_{q}}}{\left( x+n \right){{m}_{1}}!{{m}_{2}}...{{m}_{n}}!}}a_{1}^{{{m}_{1}}}a_{2}^{{{m}_{2}}}...a_{n}^{{{m}_{n}}}.$$
$A$-expansion,
$$g_{n}^{\left( m \right)}=\sum\limits_{n}{\frac{m{{\left( m+n \right)}_{q}}}{\left( m+n \right){{m}_{1}}!{{m}_{2}}...{{m}_{n}}!}}a_{1}^{{{m}_{1}}}a_{2}^{{{m}_{2}}}...a_{n}^{{{m}_{n}}},$$
is applicable to any matrix $\left( 1,xg\left( x \right) \right)$, ${{g}_{0}}=1$. It is not expansion of the binomial type, therefore we will extend the class of considered expansions. Expansions, such that  
$$g_{n}^{\left( m \right)}=\sum\limits_{n}{\frac{\left( {m}/{\varphi }\; \right){{p}_{q}}\left( \left( {m}/{\varphi }\; \right)+n \right)}{\left( \left( {m}/{\varphi }\; \right)+n \right){{m}_{1}}!{{m}_{2}}!...{{m}_{n}}!}f_{1}^{{{m}_{1}}}f_{2}^{{{m}_{2}}}...f_{n}^{{{m}_{n}}}},$$
if
$$\left[ {{x}^{n}} \right]{}_{\left( \varphi  \right)}{{A}^{m}}\left( x \right)=\sum\limits_{n}{\frac{{{p}_{q}}\left( m \right)}{{{m}_{1}}!{{m}_{2}}!...{{m}_{n}}!}f_{1}^{{{m}_{1}}}f_{2}^{{{m}_{2}}}...f_{n}^{{{m}_{n}}}},$$
where $_{\left( \varphi  \right)}A\left( x \right)$ is the generating function of the $A$-sequence of the matrix $\left( 1,x{{g}^{\varphi }}\left( x \right) \right)$, will be called the expansions of  generalized binomial type.\\
{\bfseries Theorem 7.}\emph{ $B$-expansion is the expansion of  generalized binomial type.}\\
{\bfseries Proof.} Let the matrix $\left( 1,xg\left( x \right) \right)$ is a pseudo-involution. According to the Theorem 1 and Theorem 2
$${{\left( 1,x\sqrt{g\left( x \right)} \right)}^{-1}}=\left( 1,x{{h}^{-1}}\left( x \right) \right),$$
$$h\left( x \right)=\left( 1,s\left( x \right) \right)\left( x+\sqrt{{{x}^{2}}+1} \right),  \qquad{{s}_{2n}}=0, \qquad{{s}_{2n+1}}={{{b}_{n}}}/{2}\;.$$ 
Binnomial expansion of the coefficients of the series ${{h}^{m}}\left( x \right)$ has the form
$$\left[ {{x}^{n}} \right]{{h}^{m}}\left( x \right)=\sum\limits_{n}^{{}}{\frac{{{p}_{q}}\left( m \right)}{{{m}_{0}}!{{m}_{1}}!...{{m}_{p}}!}}\frac{1}{{{2}^{q}}}b_{0}^{{{m}_{0}}}b_{1}^{{{m}_{1}}}...b_{p}^{{{m}_{p}}},$$ 
where
$$p=\left\lfloor \frac{n-1}{2} \right\rfloor , \qquad n=\sum\limits_{i=0}^{p}{{{m}_{i}}}\left( 2i+1 \right), \qquad q=\sum\limits_{i=0}^{p}{{{m}_{i}}},$$
$${{p}_{1}}\left( m \right)=m,  \qquad{{p}_{q}}\left( m \right)=m\prod\limits_{i=1}^{q-1}{\left( m+q-2i \right)}.$$ 
Corresponding expansion of the coefficients of the series ${{g}^{{m}/{2}\;}}\left( x \right)$  has the form
$$\left[ {{x}^{n}} \right]{{g}^{{m}/{2}\;}}\left( x \right)=\sum\limits_{n}^{{}}{\frac{m{{p}_{q}}\left( m+n \right)}{\left( m+n \right){{m}_{0}}!{{m}_{1}}!...{{m}_{p}}!}}\frac{1}{{{2}^{q}}}b_{0}^{{{m}_{0}}}b_{1}^{{{m}_{1}}}...b_{p}^{{{m}_{p}}}.$$
Since
$$\frac{2m}{2m+n}{{p}_{q}}\left( 2m+n \right)={{2}^{q}}m\prod\limits_{i=1}^{q-1}{\left( m+\frac{q+n}{2}-i \right)}=\frac{{{2}^{q}}m{{\left( m+k \right)}_{q}}}{m+k},$$
where 
$$k=\sum\limits_{i=0}^{p}{{{m}_{i}}\left( i+1 \right)},$$
then
$$g_{n}^{\left( m \right)}=\sum\limits_{n}^{{}}{\frac{m{{\left( m+k \right)}_{q}}}{\left( m+k \right){{m}_{0}}!{{m}_{1}}!...{{m}_{p}}!}}b_{0}^{{{m}_{0}}}b_{1}^{{{m}_{1}}}...b_{p}^{{{m}_{p}}}.$$

When deriving the $B$-expansion in section 3, we noted some its constructive properties that would be difficult to discern with a more general point of view. We note similar properties for the $A$-expansion,
$${{g}^{m}}\left( x \right)={{g}^{m-1}}\left( x \right)A\left( xg\left( x \right) \right), \qquad{{a}_{0}}=1,$$
$$g_{n}^{\left( m \right)}={{a}_{1}}\sum\limits_{i=1}^{m}{g_{n-1}^{\left( i \right)}}+{{a}_{2}}\sum\limits_{i=2}^{m+1}{g_{n-2}^{\left( i \right)}}+{{a}_{3}}\sum\limits_{i=3}^{m+2}{g_{n-3}^{\left( i \right)}+...+{{a}_{n}}}\sum\limits_{i=n}^{m+n-1}{g_{0}^{\left( i \right)}}.$$
Denote
$$\left( m|a_{1}^{{{m}_{1}}}a_{2}^{{{m}_{2}}}...a_{n}^{{{m}_{n}}} \right)=\frac{m\left( m+n-1 \right)!}{{{m}_{1}}!{{m}_{2}}...{{m}_{n}}!\left( m+n-{{m}_{1}}-{{m}_{2}}-...-{{m}_{n}} \right)!},
\quad n=\sum\nolimits_{i=1}^{n}{{{m}_{i}}}i.$$
 Then
$$\left( m|a_{r}^{{{m}_{r}}} \right)=\left[ {{x}^{{{m}_{r}}}} \right]{{\mathcal{B}}_{r}}{{\left( x \right)}^{m}}=\frac{m}{m+r{{m}_{r}}}\left( \begin{matrix}
   m+r{{m}_{r}}  \\
   {{m}_{r}}  \\
\end{matrix} \right),$$
$$\left( m|a_{1}^{{{m}_{1}}}a_{2}^{{{m}_{2}}}...a_{n}^{{{m}_{n}}} \right)=\sum\limits_{i=1}^{m}{\left( i|a_{1}^{{{m}_{1}}-1}a_{2}^{{{m}_{2}}}...a_{n}^{{{m}_{n}}} \right)}+\sum\limits_{i=2}^{m+1}{\left( i|a_{1}^{{{m}_{1}}}a_{2}^{{{m}_{2}}-1}...a_{n}^{{{m}_{n}}} \right)}+$$
$$+\sum\limits_{i=3}^{m+2}{\left( i|a_{1}^{{{m}_{1}}}a_{2}^{{{m}_{2}}}a_{3}^{{{m}_{3}}-1}...a_{n}^{{{m}_{n}}} \right)}+...+\sum\limits_{i=n}^{m+n-1}{\left( i|a_{1}^{{{m}_{1}}}a_{2}^{{{m}_{2}}}...a_{n}^{{{m}_{n}}-1} \right)},$$ 
where $\left( i|...a_{r}^{-1}... \right)=0$,
$${{g}^{m}}\left( x \right)=1+\sum\limits_{n=1}^{\infty }{\sum\limits_{n}^{{}}{\left( m|a_{1}^{{{m}_{1}}}a_{2}^{{{m}_{2}}}...a_{n}^{{{m}_{n}}} \right)}a_{1}^{{{m}_{1}}}a_{2}^{{{m}_{2}}}...a_{n}^{{{m}_{n}}}{{x}^{n}}}, \qquad n=\sum\limits_{i=1}^{n}{{{m}_{i}}}i.$$
Since
$$\left( m|a_{1}^{{{m}_{1}}}a_{2}^{{{m}_{2}}}a_{3}^{{{m}_{3}}}...a_{n}^{{{m}_{n}}} \right)=$$
$$=\frac{\left( m+n-1 \right)!m\left( m+n-{{m}_{1}}-1 \right)!}{{{m}_{1}}!\left( m+n-{{m}_{1}}-1 \right)!{{m}_{2}}!...{{m}_{n}}!\left( m+n-{{m}_{1}}-{{m}_{2}}-...-{{m}_{n}} \right)!}=$$
$$=\left( \begin{matrix}
   m+\left( n-{{m}_{1}} \right)+{{m}_{1}}-1  \\
   {{m}_{1}}  \\
\end{matrix} \right)\left( m|a_{2}^{{{m}_{2}}}a_{3}^{{{m}_{3}}}...a_{n}^{{{m}_{n}}} \right),$$
then the series ${{g}^{m}}\left( x \right)$ can also be represented in the form
$${{g}^{m}}\left( x \right)=\frac{1}{{{\left( 1-{{a}_{1}}x \right)}^{m}}}+\sum\limits_{n=1}^{\infty }{\sum\limits_{n}^{{}}{\left( m|a_{2}^{{{m}_{2}}}a_{3}^{{{m}_{3}}}...a_{n}^{{{m}_{n}}} \right)}a_{2}^{{{m}_{2}}}a_{3}^{{{m}_{3}}}...a_{n}^{{{m}_{n}}}\frac{{{x}^{n}}}{{{\left( 1-{{a}_{1}}x \right)}^{m+n}}}},
\qquad n=\sum\limits_{i=2}^{n}{{{m}_{i}}}i.$$
{\bfseries Example 6.}
$$A\left( x \right)=1+{{a}_{1}}x+{{a}_{r}}{{x}^{r}},$$  
$${{g}^{m}}\left( x \right)=\frac{1}{{{\left( 1-{{a}_{1}}x \right)}^{m}}}+\sum\limits_{n=1}^{\infty }{\left( m|a_{r}^{n} \right)a_{r}^{n}\frac{{{x}^{nr}}}{{{\left( 1-{{a}_{1}}x \right)}^{m+nr}}}}=$$
$$=\left( \frac{1}{{{\left( 1-{{a}_{1}}x \right)}^{m}}},\frac{{{a}_{r}}{{x}^{r}}}{{{\left( 1-{{a}_{1}}x \right)}^{r}}} \right){{\mathcal{B}}_{r}}{{\left( x \right)}^{m}}.$$
In particular,
$$\left( \frac{1}{{{\left( 1-{{a}_{1}}x \right)}^{m}}},\frac{{{a}_{2}}{{x}^{2}}}{{{\left( 1-{{a}_{1}}x \right)}^{2}}} \right){{\left( \frac{1-\sqrt{1-4x}}{2x} \right)}^{m}}={{\left( \frac{1-{{a}_{1}}x-\sqrt{{{\left( 1-{{a}_{1}}x \right)}^{2}}-4{{a}_{2}}{{x}^{2}}}}{2{{a}_{2}}{{x}^{2}}} \right)}^{m}}.$$

E-mail: {evgeniy\symbol{"5F}burlachenko@list.ru}
\end{document}